\documentclass{ifacconf}

\usepackage{graphicx}      
\usepackage{natbib}        
\usepackage{lipsum}
\usepackage{amsfonts}
\usepackage{graphicx}
\usepackage{epstopdf}
\usepackage[ruled,vlined,algo2e]{algorithm2e}
\usepackage{mathtools}
\usepackage{enumerate}
\usepackage{enumitem}

\DeclareMathOperator{\diag}{diag}

\begin{document}
\begin{frontmatter}

\title{Geometric statistics with subspace structure preservation for SPD matrices
\thanksref{footnoteinfo}} 

\thanks[footnoteinfo]{C.M. was supported by a Presidential Postdoctoral Fellowship at Nanyang Technological University (NTU Singapore) and an Early Career Research Fellowship at the University of Cambridge.
G.V.G. was supported by the UCL Centre for Doctoral Training in Data Intensive Science funded by STFC, and by an Overseas Research Scholarship from UCL. The research leading to these results has also received funding from the European Research Council under the Advanced ERC Grant Agreement SpikyControl n.101054323.}

\author[First]{Cyrus Mostajeran} 
\author[Second]{Nathaël Da Costa} 
\author[Third]{Graham Van Goffrier}
\author[Fourth]{Rodolphe Sepulchre}

\address[First]{School of Physical and Mathematical Sciences, Nanyang Technological University, Singapore (e-mail: cyrussam.mostajeran@ntu.edu.sg)}
\address[Second]{Tübingen AI Center, University of Tübingen, Tübingen, Germany  (e-mail: nathael.da-costa@uni-tuebingen.de)}
\address[Third]{Department of Physics and Astronomy, University College London, London, United Kingdom (e-mail: vangoffrier@gmail.com)}
\address[Fourth]{Department of Engineering, University of Cambridge, United Kingdom \& Department of Electiral Engineering, KU Leuven, Belgium (e-mail:rs771@cam.ac.uk)}

\begin{abstract}                
We present a geometric framework for the processing of SPD-valued data that preserves subspace structures and is based on the efficient computation of extreme generalized eigenvalues. This is achieved through the use of the Thompson geometry of the semidefinite cone. We explore a particular geodesic space structure in detail and establish several properties associated with it. Finally, we review a novel inductive mean of SPD matrices based on this geometry.
\end{abstract}


\end{frontmatter}

\section{Introduction}

The space of $n\times n$ symmetric positive definite (SPD) matrices  $\mathbb{S}^n_{++}$ forms the interior of the convex cone of positive semidefinite matrices in the vector space of $n\times n$ symmetric matrices. Several Riemannian geometries on SPD matrices have been proposed and used effectively in a variety of applications in computer vision, medical data analysis, machine learning, and optimization. In particular, the affine-invariant Riemannian metric---so-called because it is invariant to affine transformations of the underlying spacial coordinates---has received considerable attention in recent years and applied successfully to problems such as EEG signal processing in BCI where it has been shown to be superior to classical techniques based on feature vector classification~\citep{Barachant2012}.

The affine-invariant Riemannian metric endows the space of SPD matrices of a given dimension with the structure of a Riemannian symmetric space and a metric space of non-positive curvature~\citep{Said2017}. Computing standard geometric objects such as distances, geodesics, Riemannian exponentials and logarithms in this geometry often amounts to the computation of the generalized eigenspectrum of a pair of SPD matrices, which typically means a significant increase in computational complexity, particularly for larger matrices. In particular, the affine-invariant Riemannian structure induces the distance function $d_R:\mathbb{S}^n_{++}\times\mathbb{S}^n_{++}\rightarrow \mathbb{R}$ 
\begin{equation} \label{R distance}
d_R(X,Y)=\left(\sum_{i=1}^n\log^2\lambda_i(YX^{-1})\right)^{1/2},
\end{equation}
where $\lambda_i(YX^{-1})=\lambda_i(X^{-1/2}YX^{-1/2})$ denote the $n$ real and positive eigenvalues of $YX^{-1}$. (\ref{R distance}) satisfies two key symmetries that contribute to its success in a variety of applications: (i) affine-invariance, i.e., invariance under congruence transformations:
$d_R(X,Y)=d_R(AXA^T,AYA^T)$ for any invertible matrix $A\in\mathrm{GL}(n,\mathbb{R})$, where $A^T$ denotes the transpose of $A$, and (ii) inversion-invariance, whereby $d_R(X,Y)=d_R(X^{-1},Y^{-1})$. 

The metric (\ref{R distance}) also provides a way of defining the mean of a collection of $N$ SPD matrices $\{Y_i\}$ as the Riemannian barycentre, i.e. as
\begin{equation} \label{R mean}
    \operatorname{argmin}_{X}\sum_{i=1}^Nd_R(X,Y_i)^2,
\end{equation}
which is known to exist and be unique in this geometry~\citep{Afsari2011}. 

While the Riemannian barycentre (\ref{R mean}) has been employed effectively in numerous applications, it suffers from several drawbacks in specific contexts. Firstly, its computation can be costly for very large matrices. Secondly, it typically corrupts any sparsity or subspace structure inherent in the matrices $\{Y_i\}$, which may be an issue if such structures are of relevance to the underlying application. In this paper, we will review how a statistical framework based on the Thompson geometry of the semidefinite cone can be used to achieve affine-invariance and define interpolations and means that preserve subspace structures. These have favourable computational properties for large matrices. 
We will also establish new properties of a distinguished geodesic of the Thompson metric. In Section \ref{sec:det}, we consider the proximity of this distinguished Thompson geodesic to the affine-invariant Riemannian geodesic. We then establish in Section \ref{sec:RT distances} general bounds for the distance between their associated midpoints.

\section{Hilbert and Thompson metrics}

A subset $K$ of a vector space $V$ is called a cone if it is convex, $\mu K \subseteq K$ for all $\mu \geq 0$, and $K\cap(-K)=\{0\}$. It is said to be a closed cone if it is a closed set in $V$ with respect to the standard topology. A cone $K$ in a vector space $V$ induces a partial ordering on $V$ given by $x\leq y$ if and only if $y-x\in K$, and this partial ordering can be used to define the Hilbert projective metric and Thompson metric on $K$~\citep{Thompson1963}.

In the case of the positive semidefinite cone, we find that for strictly positive definite matrices $X, Y\succ 0$, the Hilbert and Thompson metrics take the forms
\begin{equation} \label{Hilbert matrix}
    d_H(X,Y)=\log\left(\frac{\lambda_{M}(YX^{-1})}{\lambda_{m}(YX^{-1})}\right) 
\end{equation}
and
\begin{equation} \label{Thompson matrix}
    d_T(X,Y)=\log\left(\max\{\lambda_{M}(YX^{-1}),1/\lambda_{m}(YX^{-1})\}\right),
\end{equation}
respectively, where $\lambda_{M}(A)$ and $\lambda_{m}(A)$ denote the maximum and minimum eigenvalues of the matrix $A$. Note that both of these metrics satisfy affine-invariance and inversion-invariance. Moreover, both can be computed using extreme generalized eigenvalues~\citep{Mostajeran2020,VanGoffrier2021}, which can be accessed efficiently using techniques such as Krylov subspace methods based on matrix-vector products.
Finally, we note that the Thompson metric can equivalently be expressed as
\begin{align*} 
    d_{T}(X,Y)&=\max_{1\leq i \leq n}|\log\lambda_i(YX^{-1})| \\
    &=\max\{\log\lambda_{M}(YX^{-1}),\log\lambda_{M}(XY^{-1})\}.
\end{align*}

\section{Geodesics} \label{sec:geodesics}

If $(M,d)$ is a metric space, a map $\gamma:[0,1]\rightarrow M$ is said to be a geodesic path from $x_0=\gamma(0)$ to $x_1=\gamma(1)$ if $d(\gamma(s),\gamma(t))=(t-s)d(x_0,x_1)$, whenever $0\leq s < t \leq 1$. The image of a geodesic path is called a geodesic, and a metric space is said to be a geodesic space if there exists a geodesic path joining any two points. The curve $\gamma:[0,1]\rightarrow \mathbb{S}^n_{++}$ defined by
\begin{equation} \label{Riemannian geodesic}
\gamma(t)= X\#_t Y \vcentcolon = X^{1/2}(X^{-1/2}YX^{-1/2})^t X^{1/2}
\end{equation}
is a geodesic path from $X$ to $Y$ for both the affine-invariant Riemannian metric and the Thompson metric in $\mathbb{S}^n_{++}$. In the Riemannian case, the corresponding geodesic is the unique geodesic connecting $X$ to $Y$. However, in the Thompson geometry, there is typically a family of geodesics that generally consists of an infinite number of curves connecting a pair of points in $\mathbb{S}^n_{++}$~\citep{Nussbaum1994}. 

We will use one particular geodesic in $(\mathbb{S}^n_{++},d_T)$ from $X$ to $Y$ that has particularly attractive properties. This geodesic will be denoted by $X*_t Y$ and is given by
\begin{equation} \label{Nussbaum geodesic}
X*_t Y =
\left(\frac{\lambda_M^t-\lambda_m^t}{\lambda_M-\lambda_m}\right)Y+\left(\frac{\lambda_M\lambda_m^t-\lambda_m\lambda_M^t}{\lambda_M-\lambda_m}\right)X 
\end{equation}
if $\lambda_M\neq\lambda_m$, and $X*_t Y =\lambda_m^t X$ otherwise.
Note that $X*_t Y$ reduces to a linear combination of $X$ and $Y$ with coefficients that are nonlinear functions of the extreme generalized eigenvalues of $(X,Y)$ and $t$. Furthermore, this geodesic satisfies the following joint homogeneity property for any $t\in \mathbb{R}$ and pair of positive scalars $\mu_1$ and $\mu_2$:
    \begin{equation*}
        (\mu_1X_1)*_t(\mu_2X_2)=\mu_1^{1-t}\mu_2^t(X_1*_tX_2)
    \end{equation*}
Finally, we note that if $X,Y\in \mathbb{S}^2_{++}$, then $X\#_t Y = X*_t Y$ for all $t\in [0,1]$. That is, the $*_t$ Thompson geodesic coincides with the affine-invariant Riemannian geodesic in the case of $2\times 2$ matrices.

\subsection{Subspace structure preservation}

Since $X*_t Y$ is a linear combination of $X$ and $Y$ for any $t\in\mathbb{R}$, it lies in the linear span of $X$ and $Y$, and thus preserves any subspace structure that is common to $X$ and $Y$. For example, if $X$ and $Y$ are banded, Toeplitz, or Henkel matrices, then these structures are preserved at every point along the geodesic $X*_t Y$. This is in contrast to the Riemannian geodesic (\ref{Riemannian geodesic}), which generally corrupts such structures. A further consequence of the preservation of subspace structures by the $*_t$ Thompson geodesic is that it preserves sparsity. That is, if $X$ and $Y$ are sparse SPD matrices, then $X*_t Y$ is sparse for every $t\in \mathbb{R}$. In contrast, the Riemannian geodesic $X\#_t Y$, whose construction involves computing matrix square roots, matrix products, and matrix inverses, does not preserve sparsity. Thus, the use of Riemannian interpolation to process large sparse SPD matrices may be problematic. For instance, kernel matrices in machine learning are often built as sparse matrices to facilitate the analysis of large datasets. Applying the standard affine-invariant Riemannian geometry to process such SPD matrices will typically corrupt the sparse structure, potentially resulting in intractable computations.

\subsection{Determinants along geodesics} \label{sec:det}

Here we consider the evolution of determinants along various interpolations of SPD matrices, including Euclidean, Riemannian, and Thompson geodesics. An application in which the evolution of determinants along interpolations is of practical significance is diffusion tensor imaging (DTI), which is an imaging modality used to produce non-invasive reconstructions of brain tissue connectivity. DTI is based on the assumption that the motion of water molecules in each voxel of the image is well approximated by Brownian motion. This Brownian motion is characterized by an SPD matrix, called the diffusion tensor~\citep{Basser1994}. Interpolating diffusion tensors is a basic operation in many data processing algorithms in DTI. It has been well documented that Euclidean interpolation and averaging of tensors generally results in a tensor swelling effect~\citep{Pennec2006} where the determinant of the Euclidean average of SPD matrices is often larger than the determinant of the original SPD matrices. This phenomenon is problematic as the determinant of the diffusion tensor is a measure of the dispersion of the local diffusion process, whose increase upon averaging is physically unrealistic. To overcome this difficulty, the affine-invariant and Log-Euclidean Riemannian geometries are often used in data processing tasks in DTI such as interpolation and regularization~\citep{Log-Euclidean2006,Pennec2006}. In either of these Riemannian geometries, the swelling effect is eliminated and the determinants along the geodesics between two SPD matrices with the same determinant are constant. Indeed, for $X,Y\in \mathbb{S}^n_{++}$, we have
$\det(X\#_t Y) = \det\left(X^{1/2}(X^{-1/2}YX^{-1/2})^t X^{1/2}\right) = (\det X)^{1-t}(\det Y)^t$.
In particular, 
\begin{equation*}
    \det(X\#_{1/2} Y)=\sqrt{\det(XY)}.
\end{equation*}
Thus, if $\det X = \det Y$, then $\det(X\#_t Y)=\det X = \det Y$ for all $t\in\mathbb{R}$. More generally,
we see that the logarithm of the determinant interpolates linearly along the affine-invariant Riemannian geodesics.

\begin{figure} 
    \centering
    \includegraphics[width=1.0\linewidth]{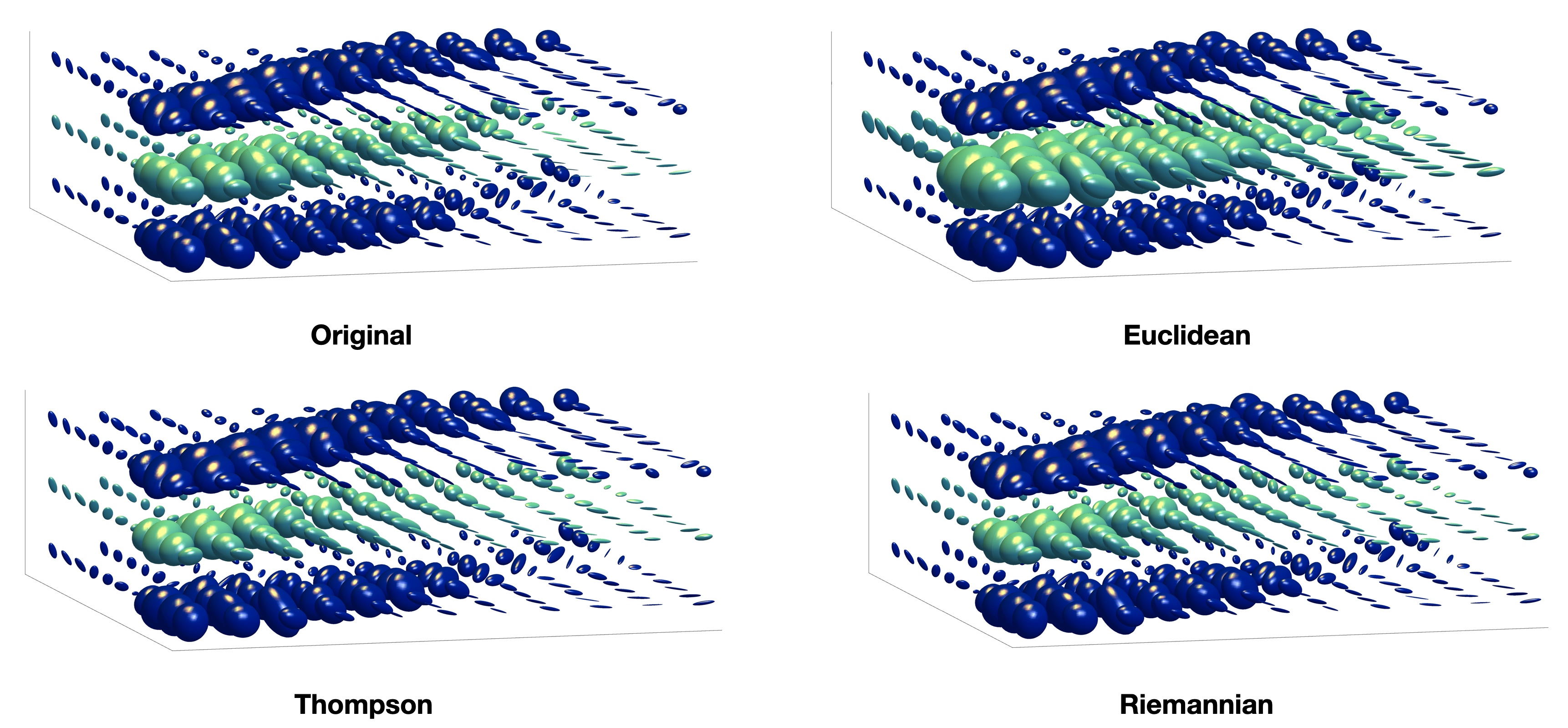}
    \caption{Reconstructions of the middle layer of an image consisting of three layers of ellipsoids representing diffusion tensors using Euclidean, Riemannian, and Thompson geodesic interpolations. The Euclidean interpolation exhibits an undesirable swelling effect. The Riemannian and Thompson interpolations generate similar results and are more faithful to the original data.}
    \label{fig:DTI}
\end{figure}

These observations raise the question of how determinants evolve along the $*_t$ Thompson geodesics. In Figure \ref{fig:DTI}, we compare interpolations using three different geodesics with real data extracted from~\cite{Camino2006} in the context of DTI. Each image depicts three layers of ellipsoids corresponding to positive definite matrices. The middle layer of the original image consisting of real data was removed and reconstructed using Euclidean, Riemannian, and Thompson interpolations of data at voxels in the top and bottom layers. We observe that in the case of $3\times 3$ matrices, the Riemannian and Thompson interpolations produce very similar reconstructions. The swelling effect in the case of the Euclidean interpolation is also clearly visible.

\begin{figure*}
    \centering
    \includegraphics[width=0.67\linewidth]{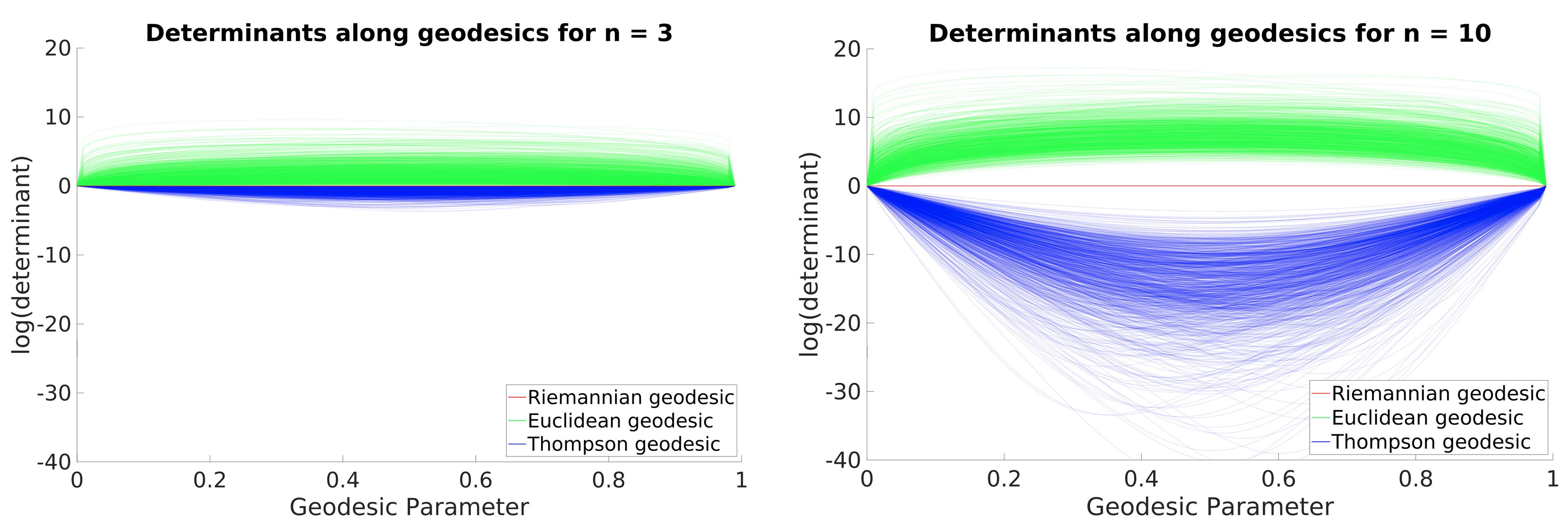}
    \caption{Plots of the logarithm of the determinants along Euclidean, Riemannian, and Thompson geodesics connecting 1,000 pairs of randomly generated $n\times n$ SPD matrices of unit determinant for $n=3$ and $n=10$. We observe that the determinants along the Riemannian geodesics are constant and equal to one as expected. The determinants along Euclidean interpolations exhibit swelling, whereas an opposite shrinkage effect is observed for interpolations along the Thompson geodesics. Both of these effects are magnified in higher dimensions.}
    \label{fig:det}
\end{figure*}

Further simulations suggest that in contrast to the Euclidean case, we seem to encounter a `shrinkage phenomenon' when using the $*_t$ Thompson geodesics, whereby determinants along the geodesic between two SPD matrices with the same determinant are reduced in value. While this effect appears modest in low dimensions, it can become quite dramatic in higher dimensions as seen in Figure \ref{fig:det}.
In the following proposition, we provide mathematical explanations for these observations and precisely characterize the circumstances under which the shrinkage phenomenon is not present along interpolations of a pair of SPD matrices.

\begin{prop}
    Let $X,Y \in \mathbb{S}^n_{++}$ and $\lambda_i$ denote the eigenvalues of $YX^{-1}$. Then
    \begin{align} \label{swelling inequality}
        &\det\left(\frac{X+Y}{2}\right) \nonumber \\
        &=\sqrt{\det(XY)}\;\prod_{i=1}^n
        \left(\frac{\sqrt{\lambda_i}+\frac{1}{\sqrt{\lambda_i}}}{2}\right) \geq \sqrt{\det(XY)}
        \end{align}
    where equality holds if and only if $X=Y$. Furthermore,
        \begin{align} \label{shrinkage inequality}
        &\det\left(X*_{1/2}Y\right) \nonumber \\
        &=
        \sqrt{\det(XY)}\;\prod_{i=1}^n
        \left[\frac{1}{\sqrt{\lambda_{M}}+\sqrt{\lambda_{m}}}\left(\sqrt{\lambda_i}+\frac{\sqrt{\lambda_{M}\lambda_{m}}}{{\sqrt{\lambda_i}}}\right)\right] \nonumber \\
        &\leq \sqrt{\det(XY)} 
    \end{align}
    where equality holds if and only if $\lambda_i\in \{\lambda_{m},\lambda_{M}\}$ for $i=1,\cdots, n$, where $\lambda_{M}$ and $\lambda_{m}$ denote the maximum and minimum eigenvalues of $YX^{-1}$, respectively.
\end{prop}
\begin{pf}
    $X^{-1/2}YX^{-1/2}$ admits a diagonalization as $QDQ^T$, where $Q$ is an orthogonal matrix and $D=\diag(\lambda_1,\cdots,\lambda_n)$. We have
    \begin{align*}
        \det\left(\frac{X+Y}{2}\right) &= 
        \det\left(X^{1/2}\frac{(I+X^{-1/2}YX^{-1/2})}{2}X^{1/2}\right) \\
        &= \det(X)\det\left(\frac{I+QDQ^T}{2}\right) \\
        &= \det(X)\det\left(\frac{I+D}{2}\right) \\
        &= \det(X)\det(D^{1/2})\det\left(\frac{D^{1/2}+D^{-1/2}}{2}\right) \\
        &= \sqrt{\det(XY)}\det\left(\frac{D^{1/2}+D^{-1/2}}{2}\right) \\
        &= \sqrt{\det(XY)}\,\prod_{i=1}^n
        \left(\frac{\sqrt{\lambda_i}+\frac{1}{\sqrt{\lambda_i}}}{2}\right),
    \end{align*}
    where we have repeatedly used the fact that $\det(AB)=\det(A)\det(B)$ for square $n\times n$ matrices $A$ and $B$. The inequality (\ref{swelling inequality}) now follows by noting the identity 
    \begin{equation*}
       \frac{1}{2}\left(\sqrt{\lambda_i}+\frac{1}{\sqrt{\lambda_i}}\right) \geq 1,
    \end{equation*}
    where equality holds if and only if $\lambda_i = 1$. Thus, each term in the product in (\ref{swelling inequality}) results in a scaling of the geometric mean $\sqrt{\det(XY)}$ by a factor greater than or equal to $1$, which contributes to the swelling effect. For the swelling effect to completely disappear, we would need $\lambda_i = 1$ for each $i$, which is only possible in the trivial case when $X=Y$.

    We now apply a similar analysis to $\det\left(X*_{1/2}Y\right)$ to find
    \begin{align*}
        &\det\left(X*_{1/2}Y\right) \\
        &=\det\left(\frac{1}{\sqrt{\lambda_{M}}+\sqrt{\lambda_{m}}}(Y+\sqrt{\lambda_{M}\lambda_{m}}X)\right) \\
        &=\det(X)\det\left(\frac{1}{\sqrt{\lambda_{M}}+\sqrt{\lambda_{m}}}(X^{-\frac{1}{2}}YX^{-\frac{1}{2}}+\sqrt{\lambda_{M}\lambda_{m}}I)\right) \\
        &=\det(X)\det\left(\frac{1}{\sqrt{\lambda_{M}}+\sqrt{\lambda_{m}}}(QDQ^T+\sqrt{\lambda_{M}\lambda_{m}}I)\right) \\
        &=\det(X)\det\left(\frac{1}{\sqrt{\lambda_{M}}+\sqrt{\lambda_{m}}}(D+\sqrt{\lambda_{M}\lambda_{m}}I)\right) \\
        &=\sqrt{\det(XY)}\det\left(\frac{1}{\sqrt{\lambda_{M}}+\sqrt{\lambda_{m}}}(D^{\frac{1}{2}}+\sqrt{\lambda_{M}\lambda_{m}}D^{-\frac{1}{2}})\right) \\
        &=\sqrt{\det(XY)}\,\prod_{i=1}^n
        \left[\frac{1}{\sqrt{\lambda_{M}}+\sqrt{\lambda_{m}}}\left(\sqrt{\lambda_i}+\frac{\sqrt{\lambda_{M}\lambda_{m}}}{{\sqrt{\lambda_i}}}\right)\right].
    \end{align*}
    Now, the inequality (\ref{shrinkage inequality}) follows from
    \begin{align*}
        \frac{1}{\sqrt{\lambda_{M}}+\sqrt{\lambda_{m}}}&\left(\sqrt{\lambda_i}+\frac{\sqrt{\lambda_{M}\lambda_{m}}}{{\sqrt{\lambda_i}}}\right) -1 \\
        &= 
        \frac{(\sqrt{\lambda_{M}}-\sqrt{\lambda_i})(\sqrt{\lambda_{m}}-\sqrt{\lambda_i})}{\sqrt{\lambda_i}(\sqrt{\lambda_{M}}+\sqrt{\lambda_{m}})} \leq 0,
    \end{align*}
    for each $i=1,\cdots, n$, which holds since $\lambda_{m} \leq \lambda_i \leq \lambda_{M}$. Moreover, we see that equality holds if and only if $\lambda_i = \lambda_{m}$ or $\lambda_i=\lambda_{M}$ for all $i$. Thus, we see that each term in the product results in a scaling of the geometric mean $\sqrt{\det(XY)}$ by a factor less than or equal to $1$, which contributes to the shrinkage effect. As the dimension increases, we expect this product to include more and more factors smaller than $1$, which explains the intensification of this shrinkage effect with the dimension $n$. The shrinkage is fully eliminated if and only if each of the factors in the product is equal to $1$, which occurs precisely if $\lambda_i\in \{\lambda_{m},\lambda_{M}\}$ for all $i=1,\cdots, n$.
\end{pf}

\subsection{Distances between geodesic mid-points}
\label{sec:RT distances}

We have seen that the affine-invariant Riemannian and $*_t$ Thompson geodesics coincide in $\mathbb{S}^2_{++}$. Moreover, simulations in $\mathbb{S}^3_{++}$ suggest that while the two geodesics no longer agree for $3\times 3$ SPD matrices, they result in qualitatively similar interpolations. This raises the question of how close these two geodesics may be in higher dimensions. To investigate this, we consider the affine-invariant Riemannian distance between the midpoints of the two geodesics connecting a pair of arbitrary SPD matrices $X$ and $Y$. To simplify the notation, we denote the Riemannian and $*_t$ Thompson geodesic midpoints by $X\#Y$ and $X*Y$, respectively.

\begin{prop}
   Let $X,Y\in\mathbb{S}^n_{++}$ and denote the eigenvalues of $YX^{-1}$ by $\lambda_i$. Then,
    \begin{align}  \label{RT distance}
        d_R&(X\#Y,X*Y) \nonumber \\
        &= 
        \left[\sum_{i=1}^n\log^2\left(\frac{1}{\sqrt{\lambda_{M}}+\sqrt{\lambda_{m}}}\left(\sqrt{\lambda_i}+\frac{\sqrt{\lambda_{M}\lambda_{m}}}{{\sqrt{\lambda_i}}}\right)\right)\right]^{1/2},
    \end{align}
    where
    $\lambda_{m}=\min_i\lambda_i$, and $\lambda_{M}=\max_i\lambda_i$. Furthermore,
    \begin{equation} \label{distance bounds}
        0 \leq d_R(X\#Y,X*Y) \leq \sqrt{n-2}\,\log\left(\cosh\frac{d_T(X,Y)}{2}\right),
    \end{equation}
    where $d_R(X\#Y,X*Y)= 0$ if and only if $\lambda_i\in\{\lambda_{m},\lambda_{M}\}$ for $i=1,\cdots,n$, and the upper bound in (\ref{distance bounds}) is attained if and only if $\lambda_{m}=1/\lambda_{M}$ and $\lambda_i = 1$ for the remaining $n-2$ eigenvalues of $YX^{-1}$.
\end{prop}
\begin{pf}
    By affine-invariance of $d_R$, we have
    \begin{align} \label{Frobenius expression}
        &d_R(X\#Y,X*Y) \nonumber \\
        &= d_R\left(X^{1/2}(X^{-1/2}YX^{-1/2})^{1/2}X^{1/2},\frac{Y+\sqrt{\lambda_{M}\lambda_{m}}X}{\sqrt{\lambda_{M}}+\sqrt{\lambda_{m}}}\right) \nonumber \\ &= d_R\left((X^{-1/2}YX^{-1/2})^{1/2},\frac{X^{-1/2}YX^{-1/2}+\sqrt{\lambda_{M}\lambda_{m}}\,I}{\sqrt{\lambda_{M}}+\sqrt{\lambda_{m}}}\right) \nonumber \\
        &= \bigg\|\log\left(\frac{1}{\sqrt{\lambda_{M}}+\sqrt{\lambda_{m}}}\,\left(\Sigma^{1/2}+\sqrt{\lambda_{M}\lambda_{m}}\Sigma^{-1/2}\right)\right)\bigg\|,
    \end{align}
    where $\Sigma=X^{-1/2}YX^{-1/2}$ and $\|\cdot\|$ denotes the Frobenius norm. By noting that $\Sigma$ admits a diagonalization $\Sigma=QDQ^T$, where $Q$ is an orthogonal matrix and $D=\diag(\lambda_1,\cdots,\lambda_n)$, and considering the eigenvalues of the matrix logarithm in (\ref{Frobenius expression}), we arrive at (\ref{RT distance}).

    It follows from (\ref{RT distance}) that $d_R(X\#Y,X*Y)=0$ if and only if 
    \begin{equation*}
        \frac{1}{\sqrt{\lambda_{M}}+\sqrt{\lambda_{m}}}\left(\sqrt{\lambda_i}+\frac{\sqrt{\lambda_{M}\lambda_{m}}}{{\sqrt{\lambda_i}}}\right) = 1
    \end{equation*}
    for $i=1,\cdots,n$, which holds if and only if $(\sqrt{\lambda_i}-\sqrt{\lambda_{M}})(\sqrt{\lambda_i}-\sqrt{\lambda_{m}})=0$. That is, $d_R(X\#Y,X*Y)= 0$ if and only if $\lambda_i\in\{\lambda_{m},\lambda_{M}\}$ for $i=1,\cdots,n$. 
    
    Now denote by $r$ the Thompson distance $d_T(X,Y)=\max_i|\log\lambda_i|$ between $X$ and $Y$. It follows that $\lambda_{M}=e^r$ or $\lambda_{m}=e^{-r}$. Suppose that $\lambda_{M}=e^r$ and note that (\ref{RT distance}) can be expressed as
    \begin{equation} \label{rho sum}
        d_R(X\#Y,X*Y)=\left(\sum_{i=1}^n \rho(\lambda_i)\right)^{1/2},
    \end{equation}
    where 
    \begin{equation*}
        \rho(x)=\log^2\left(\frac{1}{\sqrt{\lambda_{M}}+\sqrt{\lambda_{m}}}\left(\sqrt{x}+\frac{\sqrt{\lambda_{M}\lambda_{m}}}{{\sqrt{x}}}\right)\right),
    \end{equation*}
    for $x\in [\lambda_{m},\lambda_{M}]$. To establish the upper bound in (\ref{distance bounds}), we seek to maximize (\ref{RT distance}) subject to the constraint that $e^{-r}\leq \lambda_{m}\leq \lambda_i \leq \lambda_{M}=e^r$. For a given $\lambda_{m}$, this is achieved by maximizing $\rho(\lambda_i)$ for every choice $\lambda_i$ not corresponding to the extreme eigenvalues. $\rho$ is a non-negative continuous function on $[\lambda_{m},\lambda_{M}]$ that is smooth on $(\lambda_{m},\lambda_{M})$ and satisfies $\rho'(x)=0$ if and only if $x=\sqrt{\lambda_{m}\lambda_{M}}$. Since $\rho(\lambda_{m})=\rho(\lambda_{M})=0$, it follows that (\ref{rho sum}) is maximized for a given $\lambda_{m}$ and $\lambda_{M}$ if and only if all remaining eigenvalues satisfy $\lambda_i=\sqrt{\lambda_{m}\lambda_{M}}$. For such a distribution of eigenvalues, (\ref{rho sum}) reduces to
    \begin{equation*}
        d_R(X\#Y,X*Y)=\left((n-2)\log^2\left(\frac{2e^{r/4}\lambda_{m}^{1/4}}{e^{r/2}+\lambda_{m}^{1/2}}\right)\right)^{1/2},
    \end{equation*}
    which in turn is maximized as a function of $\lambda_{m}$ subject to the constraint $e^{-r}\leq \lambda_{m} \leq \lambda_{M}=e^r$ when $\lambda_{m}=e^{-r}$, at which point
    \begin{align*}
        d_R(X\#Y,X*Y)&=\left((n-2)\log^2\left(\frac{2}{e^{r/2}+e^{-r/2}}\right)\right)^{1/2} \\
        &= \sqrt{n-2}\,\log\left(\cosh\frac{r}{2}\right).
    \end{align*}
    Repeating the analysis by starting with the assumption that $\lambda_{m}=e^{-r}$ instead of $\lambda_{M}=e^{r}$ would yield the same conclusion. That is, the upper bound in (\ref{distance bounds}) holds and is attained if and only if $\lambda_{m}=e^{-r}$, $\lambda_{M}=e^r$, and $\lambda_i = \sqrt{\lambda_{m}\lambda_{M}}=1$ for the remaining $n-2$ eigenvalues of $YX^{-1}$, where $r=d_T(X,Y)$.
\end{pf}

\begin{figure*}[t]
    \centering
    \includegraphics[width=0.75
\linewidth]{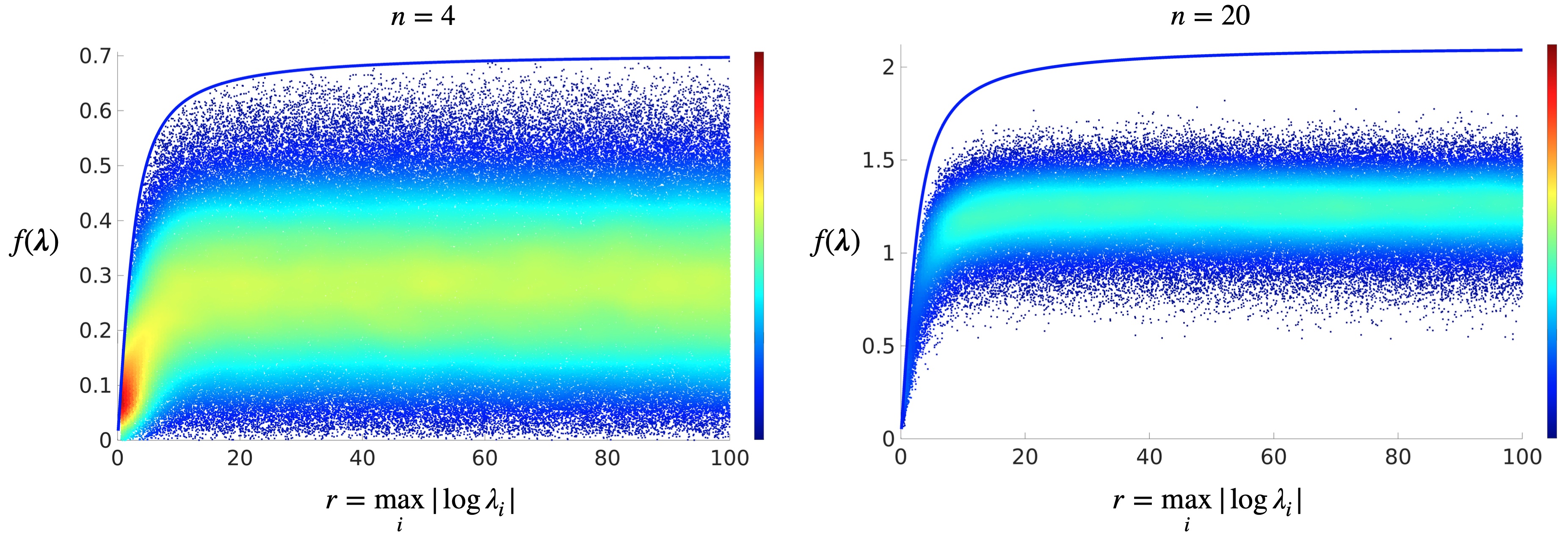}
    \caption{Plots of $f(\boldsymbol{\lambda})$ (\ref{f}) against $r=\max_i|\log\lambda_i|$ for 100,000 samples $\boldsymbol{\lambda}\in\mathbb{R}^n_+$ for $n=4$ and $n=20$. The solid blue curves represent the upper bounds of $f(\boldsymbol{\lambda})$ given in (\ref{f bounds}). The color scheme is used to indicate the density of plot points $(r,f(\boldsymbol{\lambda}))$.}
    \label{fig:RT distance plots}
\end{figure*}

While we see that the bounds in (\ref{distance bounds}) are attainable, it is generally unlikely that two SPD matrices $X$ and $Y$ will have the generalized eigenvalues required to attain these bounds. Indeed, the required distributions of generalized eigenvalues become increasingly unlikely as the dimension $n$ grows. To develop a more practical sense of the typical values of the distance $d_R(X\#Y,X*Y)$, we can compute this distance for a large number of randomly generated SPD matrices. As (\ref{RT distance}) provides an expression for $d_R(X\#Y,X*Y)$ in terms of the eigenvalues of $YX^{-1}$, it allows us to efficiently simulate a large number of matrix distance computations by sampling vectors of generalized eigenvalues instead of generating a large number of pairs of SPD matrices and computing the corresponding generalized eigenvalues. Thus, we consider the `normalized' distance function 
\begin{align} \label{f}
    f(\boldsymbol{\lambda}) \vcentcolon=& \frac{\sqrt{\sum_{i}\log^2\left(\frac{1}{\sqrt{\lambda_{M}}+\sqrt{\lambda_{m}}}\left(\sqrt{\lambda_i}+\frac{\sqrt{\lambda_{M}\lambda_{m}}} {{\sqrt{\lambda_i}}}\right)\right)}}{\max_i|\log\lambda_i|} \nonumber \\
    =&\frac{d_R(X\#Y,X*Y)}{d_T(X,Y)}
\end{align}
where $f$ is a positive-valued function of $\boldsymbol{\lambda}=(\lambda_1,\cdots,\lambda_n)$. The results in (\ref{fig:RT distance plots})
were produced by generating 100,000 points $\boldsymbol{\lambda}$ for values of $r=\max_i|\log\lambda_i|$ from 0 to 100 and computing the corresponding $f(\boldsymbol{\lambda})$ in the cases $n=4$ and $n=20$. The solid blue curve above each plot indicates the corresponding upper bound as a function of $r$. We note that in either case most points concentrate away from the bounds and increasingly so in higher dimensions. By (\ref{distance bounds}), we have 
\begin{equation} \label{f bounds}
    0 \leq f(\boldsymbol{\lambda})\leq \sqrt{n-2}\,\frac{\log\left(\cosh\frac{r}{2}\right)}{r}.
\end{equation}

\section{Inductive mean}

An important step in developing a computational framework for performing analysis and statistics on SPD-valued data using extreme generalized eigenvalues is to provide a suitable definition for the mean of a collection of $k$ SPD matrices whose computation can be based primarily on finding a sequence of extreme generalized eigenvalues. Here, we review a result from~\cite{Mostajeran2023} that provides an inductive algorithm based on Thompson geodesics that achieves this.

\begin{algorithm2e}[h]
\SetAlgoLined
\For{$i\geq 1$} {
Set $j \equiv i \mod k$ for $1\leq j \leq k$. \label{line3} \\
Define $X_{i+1} = X_i*_{\frac{1}{i+1}}Y_j$.
}
\Return{$(X_1,X_2,X_3,\cdots)$.}
\caption{Generate inductive sequence of SPD matrices $(X_i)_{i\geq 1}$ from an initial point $X_1$ and the finite ordered set $\mathcal{P}=(Y_1,\cdots,Y_k)\subset\mathbb{S}^n_{++}$}
\label{alg:inductive mean}
\end{algorithm2e}

\begin{thm}\label{cvgence_thm} Let $(X_i)_{i\geq 1}$ denote any sequence generated by Algorithm \ref{alg:inductive mean}. Then, $(X_i)_{i\geq 1}$ converges to a unique point $X^*$ that is independent of the choice of initialization $X_1$ and lies in the linear span of $\{Y_1,\cdots,Y_k\}$.
\end{thm}

Thanks to Theorem \ref{cvgence_thm}, we have a well-defined inductive mean $M(Y_1,\dots,Y_k)= X^*$ of $\mathcal{P}=(Y_1,\cdots,Y_k)\subset\mathbb{S}^n_{++}$. This mean is invariant to permutation of the $Y_j$ and satisfies the following properties:
\begin{enumerate}
\item Affine-equivariance: 
\begin{equation*}
 M(AY_1A^T,\dots, AY_kA^T) = AM(Y_1,\dots,Y_k)A^T   
\end{equation*} 
for any invertible matrix $A$.
\item Joint homogeneity: 
\begin{equation*}
M(c_1Y_1,\dots,c_kY_k) = (c_1\dots\cdot c_k)^{1/k}M(Y_1,\dots,Y_k)     
\end{equation*}
for any $c_1,\dots,c_k>0$.
\end{enumerate}

See~\cite{Mostajeran2023} for the proof of Theorem \ref{cvgence_thm} and the above properties of the resulting inductive Thompson mean.

The subspace structure preservation property of the inductive mean leads to the following corollaries concerning sparsity preservation.

\begin{cor}[Sparsity preservation I]\label{sparsity_preservation_1}
If $\{Y_1\dots,Y_k\}$ is a set of SPD matrices with the same sparsity pattern (i.e., with non-zero elements restricted to a common set of entries), then $M(Y_1,\dots,Y_k)$ has the same sparsity pattern.
\end{cor}

\begin{cor}[Sparsity preservation II]\label{sparsity_preservation_2}
If $\{Y_1\dots,Y_k\}\subset \mathbb{S}_{++}^n$ is a set of sparse SPD matrices with $k<<n^2$, then $M(Y_1,\dots,Y_k)$ is sparse.
\end{cor}

\begin{figure}
    \centering
    \includegraphics[width=
0.96\linewidth]{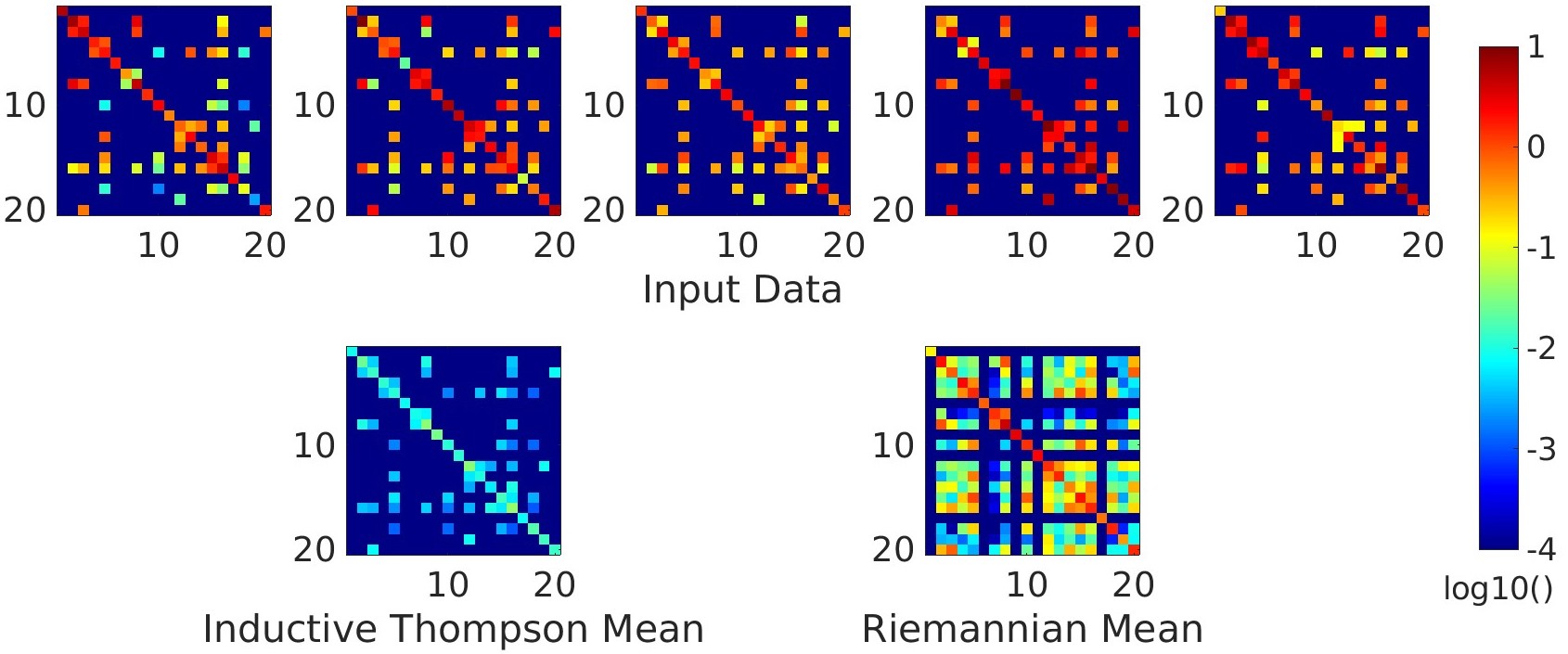}
    \caption{The inductive Thompson mean and Riemannian mean for a set of 5 input data SPD matrices with the same sparsity pattern.}
    \label{fig:inductive1}
\end{figure}

\begin{figure}
    \centering
    \includegraphics[width=
0.96\linewidth]{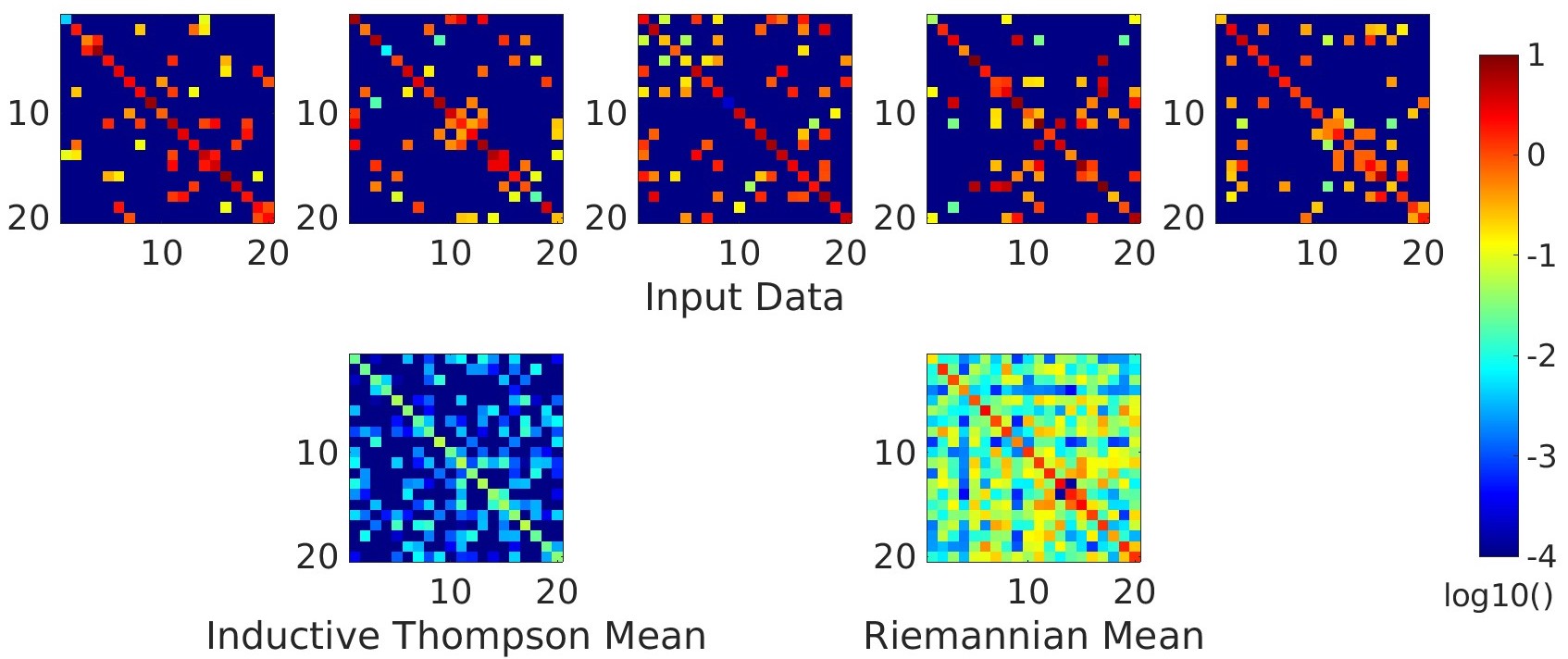}
    \caption{The inductive Thompson mean and Riemannian mean for a set of 5 input data SPD matrices with distinct sparsity patterns.}
    \label{fig:inductive2}
\end{figure}

Figure \ref{fig:inductive1} shows the inductive Thompson mean and Riemannian mean of 5 input data matrices with the same sparsity pattern, showing how the inductive Thompson mean preserves the sparsity pattern whereas the Riemannian mean corrupts it. Figure \ref{fig:inductive2} shows the results of the corresponding simulations for a second set of 5 sparse input matrices with distinct sparsity patterns.

\bibliography{ifacconf}             
                                                   







\end{document}